\documentclass[11pt,a4paper]{article}
\usepackage[margin=25mm]{geometry}
\usepackage[T1]{fontenc}
\usepackage{lmodern,microtype,amsmath,amssymb,amsthm,mathtools,mathrsfs}
\usepackage{tikz,float}
\usepackage[hidelinks]{hyperref}
\newtheorem{theorem}{Theorem}[section]
\newtheorem{lemma}[theorem]{Lemma}
\newtheorem{proposition}[theorem]{Proposition}

\theoremstyle{remark}
\numberwithin{equation}{section}
\newcommand{\R}{\mathbb R}
\newcommand{\Diss}{\mathscr D}
\newcommand{\sech}{\operatorname{sech}}
\hypersetup{pdftitle={Critical radial Green potential heat flow: a tension-controlled residue and a two-bubble example},pdfauthor={Zhe Zhang}}
\title{Critical radial Green potential heat flow:\\ a tension-controlled residue 
and a two-bubble example}
\author{Zhe Zhang}
\date{November 2026}
\begin{document}
\maketitle

\begin{abstract}
We compare two heat flows for the same corotational map of a disk
into a warped product of two spheres, with the same fixed boundary
values. The ordinary harmonic map heat flow develops a singularity
in finite time. When the domain metric is instead determined by the
critically normalized Green potential of the current energy
distribution, the solution exists for all time. Two harmonic spheres
of energy $4\pi$ form at separated scales, in addition to a
nonconstant harmonic root. The energy and the image diameter of the
necks tend to zero without restriction to a sequence of times.
The main estimate controls the difference between the inner source
mass and its critical value by the weighted tension. This is necessary
because that difference is integrated over a neck of unbounded
logarithmic length. Along a sequence of times the Green metric on
the middle neck is asymptotically cylindrical after a spatially
constant normalization.
\end{abstract}

\section{Introduction}\label{sec:intro}

Topping \cite{TwoScales} constructed a harmonic map heat flow in which
two bubbles develop at the same point at infinite time and at different
scales. His example uses two corotational components and a warped
metric on $S^2\times S^2$. That construction is the inspiration for
the example in this paper.

We modify the warping and the boundary values in order to compare two
flows from the \emph{same} initial map. Under the ordinary harmonic map
heat flow, the map develops a singularity in finite time. Under the
Green--potential heat flow, with its coupling fixed at the critical
normalization, it stays smooth at every finite time and develops two
harmonic spheres at separated scales as $t\to\infty$. In the latter
case the two spheres and a nonconstant root account for all the energy,
and the necks have vanishing image diameter. The precise statement is
Theorem~\ref{thm:example}.

There are two changes to Topping's model. The warping is constant near
both poles of the first sphere, and the second angle has a small
boundary overshoot beyond $\pi$. The overshoot permits a scalar
comparison argument for finite-time blow-up, in the spirit of
Chang, Ding and Ye \cite{CDY}. Constancy of the warping near the poles
lets us identify the limiting components and exclude a stationary
configuration in which both transitions remain at one scale.

For the modified flow, the main difficulty is the coefficient between
the two concentrating scales. The critical source mass is $2\pi$,
but a small error in that mass can accumulate over a long logarithmic
neck. We control this error by the current weighted tension. This
first prevents finite-time collapse of the scales and then supplies
one sequence of times on which both components are almost harmonic.
An energy estimate subsequently gives the two-bubble and no-neck
limits for the full time variable. The asymptotically cylindrical
Green metric is obtained only along selected times, after one spatially
constant normalization on each slice.

This example was announced at the end of Section~5 of \cite{GPHF}.
The present paper gives its proof directly; the general global
existence and bubble-tree theorems of that paper are not used.

\section{Construction of the example}\label{sec:setup}

Let $D=B_1\subset\R^2$. Choose a fixed smooth function $g:[0,\pi]\to[0,1]$
such that
\[
 g=0\ \hbox{on }[0,\pi/2],\qquad
 g=1\ \hbox{on }[3\pi/4,\pi],\qquad
 g'>0\ \hbox{on }(\pi/2,3\pi/4).
\]
Put $\mathfrak f=1+\epsilon_w g$, where $\epsilon_w>0$ is fixed and
sufficiently small, and write $f_\pi=1+\epsilon_w$.
The target is $N=S^2\times S^2$ with the smooth metric
\[
 h_N=h_{S^2_\alpha}+\mathfrak f(\alpha)h_{S^2_\beta}.
\]
Here $\alpha$ is colatitude in the first sphere; the warping is constant near
its poles. With $n(v,\theta)=(\sin v\cos\theta,\sin v\sin\theta,\cos v)$,
we write $u=(n(\alpha,\theta),n(\beta,\theta))$. The second angle is
a real lift and may pass through $\pi$. The energy and tension without
a metric subscript are for the Euclidean disk and the metric $h_N$.

\paragraph{Initial and boundary data.}
Fix $b=\pi+\delta$, with $0<\delta<\pi/4$, and impose
\[
 \alpha(0,t)=\beta(0,t)=0,\qquad
 \alpha(1,t)=\pi,\quad \beta(1,t)=b.
\]
Write
\[
 H_b(r)=\pi+2\arctan(c_b r),\qquad
 c_b=\tan(\delta/2),\qquad E_b=2\pi(1+\cos b).
\]
Choose a smooth nondecreasing $\vartheta$ equal to zero on $[0,1/2]$
and one on $[3/4,1]$. For $0<\lambda_0\ll R_0\ll1$, take
\begin{equation}
 \begin{split}
 \alpha_0(r)&=(1-\vartheta(r))\,2\arctan(r/R_0)
                         +\vartheta(r)\pi,\\
 \beta_0(r)&=H_b(r)-(1-\vartheta(r))
                         \bigl(\pi-2\arctan(r/\lambda_0)\bigr).
 \end{split}
 \label{eq:data}
\end{equation}
These are smooth equivariant maps; near the boundary they agree with
the stationary pair $(\pi,H_b)$ and satisfy every compatibility
condition. They obey $0\le\alpha_0\le\beta_0\le b$, $\alpha_0\le\pi$.

\paragraph{The Green--potential equation.}
Fix a smooth radial function $0\leq\chi\leq1$, equal to one on $B_1$
and zero outside $B_2$, and put $\chi_{p,\ell}(y)=\chi((y-p)/\ell)$.
For a nonnegative measure $\mu$ on the disk (or on the plane when
calibrating the source), put
\[
 M=\int\chi_{p,\ell}\,d\mu,\quad b_{p,\ell}=M^{-1}\int y\chi_{p,\ell}\,d\mu,
 \quad \sigma^2=M^{-1}\int|y-b_{p,\ell}|^2\chi_{p,\ell}\,d\mu .
\]
Choose smooth nondecreasing cutoffs into $[0,1]$: $\Gamma_1$ is zero
on $[0,\varepsilon_{\rm pkt}/4]$ and one on
$[\varepsilon_{\rm pkt}/2,\infty)$; $\Gamma_2$ is zero on $[0,\rho_-]$
and one on $[\rho_+,\infty)$, where $0<\rho_-<\rho_+$.
Write $W_\mu=\Gamma_1(M)\Gamma_2(\sigma/\ell)$, assigning zero
weight when $M=0$, without evaluating the center or variance.

Write $Q_\ell(r)=2\arctan(r/\ell)$ and
$Q(x)=2\arctan(e^{-x})$, so $Q_\ell(e^{-s})=Q(s+\log\ell)$.
The energy measure of $n(Q_1(r),\theta)$
is $d\mu_Q=4(1+|y|^2)^{-2}\,dy$. Define
\[
 A_Q^{\rm raw}:=\int_{\R^2\times(0,\infty)}
                  W_{\mu_Q}(p,\ell)\,\frac{dp\,d\ell}{\ell^3}.
\]
Fix
$\varepsilon_{\rm pkt}<32\pi/1681$ and $\rho_+<25/(82\sqrt2)$.
These choices give $0<A_Q^{\rm raw}<\infty$, as verified below
at the beginning of Subsection~\ref{sec:source-bands}.

Use $\mu_u=\tfrac12|du|^2dx$, $W_u=W_{\mu_u}$, a fixed
$0<r_*<1/16$, and the Dirichlet heat kernel $H^D$. Set
\[
 \nu_u(x)=\frac{2\pi}{A_Q^{\rm raw}}
 \int_{D\times(0,r_*)}W_u(p,\ell)H^D_{\ell^2}(x,b_{p,\ell})
                                      \frac{dp\,d\ell}{\ell^3}.
\]
The normalized plane parameter mass is $A_Q=2\pi$.
There is no subtraction of a spatial mean on the disk and no
renormalization by the source mass at a later time. With $\kappa=1$,
the critical exponent is $\gamma_Q=\kappa A_Q/\pi=2$. The two flows are
$u_t=\tau(u)$ and
\begin{equation}
 -\Delta\phi_u=\nu_u,\quad \phi_u|_{\partial D}=0,\qquad
 a_u=e^{-2\phi_u},\quad u_t=a_u\tau(u).
 \label{eq:critical-pde}
\end{equation}
The equations and the data respect simultaneous rotations and reflections
of the domain and both spheres. In particular, rotational equivariance makes $\nu_u$ and $\phi_u$ radial.

\begin{theorem}[Same-data comparison at critical coupling]
\label{thm:example}
Choose $\epsilon_w>0$ sufficiently small. For every sufficiently small
$\delta>0$, choose $R_0>0$ sufficiently small and then $\lambda_0>0$
sufficiently small relative to $R_0$ and $\delta$. The following hold
for the same initial and boundary data, with initial map \eqref{eq:data}.
\begin{enumerate}
\item The ordinary coupled flow develops a finite-time singularity.
\item The critical PW flow \eqref{eq:critical-pde} is smooth at every
finite time. There are smooth modulation scales
$R(t)=e^{-S_\alpha(t)}$, $\lambda(t)=e^{-S_\beta(t)}$ such that
\[
 S_\alpha(t)\longrightarrow\infty,\qquad
 L(t):=S_\beta(t)-S_\alpha(t)\longrightarrow\infty.
\]
With $s=-\log r$ and $Q(x)=2\arctan(e^{-x})$, the root and vertex
limits, locally in energy and uniformly on compact sets, are
\[
 u(t)\longrightarrow(\pi,H_b)\quad\hbox{on }D\setminus\{0\},
\]
\[
 (\alpha,\beta)(S_\alpha(t)+x,t)\longrightarrow(Q(x),\pi),
 \qquad
 (\alpha,\beta)(S_\beta(t)+x,t)\longrightarrow(0,Q(x)).
\]
These limits hold for the full time variable, and
\begin{equation}
 \lim_{t\to\infty}E(u(t))=8\pi+f_\pi E_b.
 \label{eq:energy-limit}
\end{equation}
\item Set
\[
 I_J^0=[J,S_\alpha-J],\quad
 I_J^1=[S_\alpha+J,S_\beta-J],\quad
 I_J^2=[S_\beta+J,\infty),
\]
omitting empty intervals. Then
\[
 \lim_{J\to\infty}\limsup_{t\to\infty}
 \sum_{k=0}^2 E(u(t);I_J^k\times S^1)=0,
\]
\[
 \lim_{J\to\infty}\limsup_{t\to\infty}
 \max_{0\le k\le2}\operatorname{diam}_N u(I_J^k\times S^1,t)=0.
\]
\item Along a sequence $t_j\to\infty$ the middle neck becomes cylindrical
after one spatially constant metric normalization per slice. More precisely, put
$q(t,r)=a(t,r)/r^2$ and $r_m(t)=\sqrt{\lambda(t)R(t)}$. Then
\[
 \lim_{K\to\infty}\limsup_{j\to\infty}
 \sup_{K\lambda(t_j)\le r\le R(t_j)/K}
 \left|\log\frac{q(t_j,r)}{q(t_j,r_m(t_j))}\right|=0.
\]
\end{enumerate}
\end{theorem}

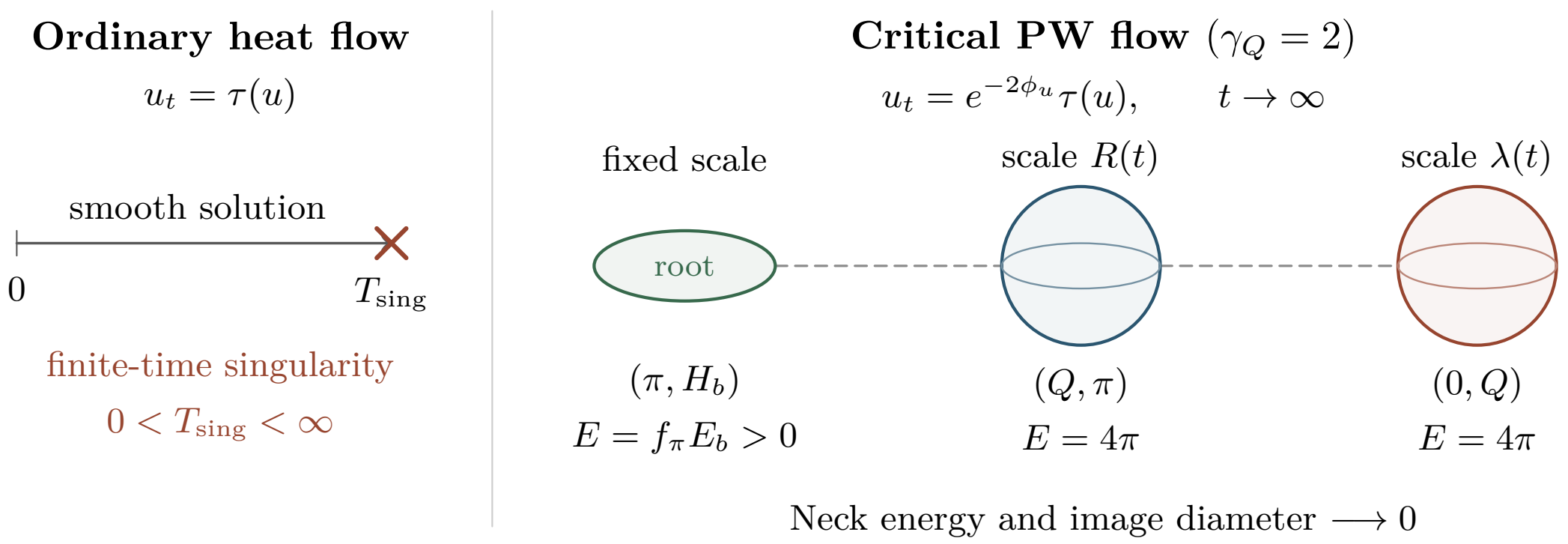
\begin{figure}[H]
\centering
\begin{tikzpicture}[x=1cm,y=1cm,font=\footnotesize,
  line cap=round,line join=round]
  \definecolor{pwblue}{RGB}{43,86,112}
  \definecolor{pwrust}{RGB}{151,69,47}
  \definecolor{pwgreen}{RGB}{55,105,76}
  \path[use as bounding box] (0,-.3) rectangle (15.2,4.7);
  \draw[black!20] (4.55,.1)--(4.55,4.65);

  % The ordinary flow: only the maximal smooth lifetime is depicted.
  \node[font=\small\bfseries] at (2.15,4.4) {Ordinary heat flow};
  \node at (2.15,3.9) {$u_t=\tau(u)$};
  \draw[black!65,line width=.65pt,->] (.35,2.6)--(3.65,2.6);
  \draw[black!65] (.35,2.49)--(.35,2.71);
  \draw[pwrust,line width=1.15pt]
    (3.53,2.47)--(3.79,2.73) (3.53,2.73)--(3.79,2.47);
  \node[above] at (1.95,2.68) {smooth solution};
  \node[below] at (.35,2.43) {$0$};
  \node[below] at (3.66,2.43) {$T_{\mathrm{sing}}$};
  \node[text=pwrust,align=center] at (2.15,1.25)
    {finite-time singularity\\[3pt]$0<T_{\mathrm{sing}}<\infty$};

  % Critical PW: rescaled limiting components, not a metric drawing.
  \node[font=\small\bfseries] at (9.95,4.4) {Critical PW flow $(\gamma_Q=2)$};
  \node at (9.95,3.9)
    {$u_t=e^{-2\phi_u}\tau(u),\qquad t\to\infty$};

  \node at (6.25,3.35) {fixed scale};
  \node at (9.75,3.35) {scale $R(t)$};
  \node at (13.25,3.35) {scale $\lambda(t)$};

  \draw[black!45,densely dashed,line width=.6pt]
    (7.05,2.4)--(9.05,2.4);
  \draw[black!45,densely dashed,line width=.6pt]
    (10.45,2.4)--(12.55,2.4);

  % The disk root is drawn without a spherical equator.
  \draw[pwgreen,fill=pwgreen!7,line width=.8pt]
    (6.25,2.4) ellipse[x radius=.8,y radius=.31];
  \node[text=pwgreen,font=\scriptsize] at (6.25,2.4) {root};

  % The two spherical domains are shown at the same schematic size.
  \draw[pwblue,fill=pwblue!6,line width=.8pt]
    (9.75,2.4) circle[radius=.7];
  \draw[pwblue!65,line width=.45pt]
    (9.75,2.4) ellipse[x radius=.7,y radius=.20];
  \draw[pwrust,fill=pwrust!6,line width=.8pt]
    (13.25,2.4) circle[radius=.7];
  \draw[pwrust!65,line width=.45pt]
    (13.25,2.4) ellipse[x radius=.7,y radius=.20];

  \node[align=center] at (6.25,1.14)
    {$(\pi,H_b)$\\[3pt]$E=f_\pi E_b>0$};
  \node[align=center] at (9.75,1.14)
    {$(Q,\pi)$\\[3pt]$E=4\pi$};
  \node[align=center] at (13.25,1.14)
    {$(0,Q)$\\[3pt]$E=4\pi$};
  \node at (9.95,.15) {Neck energy and image diameter $\longrightarrow0$};
\end{tikzpicture}
\caption{The same initial map under the two flows. On the right,
the root and the two rescaled spheres are drawn at comparable sizes,
although $R(t)\to0$ and $\lambda(t)/R(t)\to0$.
The dashed segments represent transition regions with vanishing energy
and image diameter. The normalized metric limit on the middle neck
holds along the selected times of Theorem~\ref{thm:example}.
The schematic was prepared in TikZ with AI-assisted code.}
\label{fig:same-data-comparison}
\end{figure}

Figure~\ref{fig:same-data-comparison} summarizes the comparison.
Throughout, the boundary values are fixed in time. The metric statement
in the theorem is equivalent to uniform convergence of the coefficient
of $q(t_j,r_m(t_j))\,e^{2\phi_u}(dr^2+r^2d\theta^2)$ to that of
$ds^2+d\theta^2$ on the indicated trimmed necks, in the displayed
order of limits. No convergence of metric derivatives is asserted.

\section{Proof of Theorem~\ref{thm:example}}\label{sec:proof}

We first write the two angle equations and establish their common
range bounds. We then compare the two evolutions: a scalar barrier
gives a finite-time singularity for the ordinary flow, whereas an
estimate for the actual potential source prevents it for the modified flow.
It remains to identify the only possible limiting components and
show that their necks carry no energy.

For $s=-\log r$, subscripts $s,t$ denote derivatives of the lifted
angles at fixed $t,s$, respectively. The Dirichlet energy becomes
\[
 E=\pi\int_0^\infty[\alpha_s^2+\sin^2\alpha+
                 \mathfrak f(\alpha)(\beta_s^2+\sin^2\beta)]\,ds .
\]
The fixed boundary values eliminate the boundary term in its first
variation. For fixed $t$ abbreviate the cylindrical coefficient by
$q(s)=q(t,e^{-s})=e^{2s}a_u(t,e^{-s})$ for the Green-potential
flow; for the ordinary flow $q(s)=e^{2s}$. The scalar profile
satisfies $Q'=-\sin Q=-\sech x$.

The cylindrical equations and dissipation are
\begin{align*}
 \alpha_t&=qF_\alpha,&
 F_\alpha&=\alpha_{ss}-\sin\alpha\cos\alpha
 -\tfrac12\mathfrak f'(\alpha)(\beta_s^2+\sin^2\beta),\\
 \beta_t&=qF_\beta,&
 F_\beta&=\beta_{ss}-\sin\beta\cos\beta
 +\frac{\mathfrak f'(\alpha)}{\mathfrak f(\alpha)}\alpha_s\beta_s,\\
 E'(t)&=-\Diss(t),&
 \Diss&=2\pi\int_0^\infty q
       (F_\alpha^2+\mathfrak f(\alpha)F_\beta^2)\,ds.
\end{align*}
For the following comparisons, the coefficient of a given smooth solution
is held fixed. The scalar maximum principle, applied to this positive
coefficient, preserves $0\le\alpha\le\pi$ and $0\le\beta\le b$.
At a first contact of $\beta-\alpha$ in $[0,\pi]$, its derivative
vanishes and the additional forcing in its equation is
\[
 \mathfrak f'(\alpha)
 \left[\left(\frac1{\mathfrak f(\alpha)}+\frac12\right)\alpha_s^2
                         +\frac12\sin^2\alpha\right]\ge0.
\]
Consequently $\alpha\le\beta$ is preserved. No comparison of two
different nonlocal coefficients is involved.

\subsection{The ordinary flow develops a finite-time singularity}\label{sec:ordinary}

We first prove the ordinary-flow assertion of Theorem~\ref{thm:example}. In physical radius
the ordinary equations are
\begin{align*}
 \alpha_t&=\mathcal T\alpha-\tfrac12\mathfrak f'(\alpha)
                           (\beta_r^2+r^{-2}\sin^2\beta),\\
 \beta_t&=\mathcal T\beta+
                \frac{\mathfrak f'(\alpha)}{\mathfrak f(\alpha)}\alpha_r\beta_r,
 \qquad
 \mathcal T v=v_{rr}+r^{-1}v_r-r^{-2}\sin v\cos v.
\end{align*}
All comparisons below are restricted to the common smooth lifetime. The
range bounds remain valid. Since $\mathfrak f'\ge0$, the first
angle is a scalar subsolution. Choose
\[
 K_1\ge\max\{\pi,\sup_{r>0}\alpha_0(r)/r\},
 \qquad K'=\tfrac23K^3,\quad K(0)=K_1.
\]
The function $K(t)r$ is a scalar supersolution on the smooth lifetime because
\[
 x-\sin x\cos x=2\int_0^x\sin^2y\,dy\le\tfrac23x^3.
\]
Thus $\alpha\le2K_1r$ for
$t\le\tau_\alpha=3/(8K_1^2)$. With $r_c=1/(16K_1)$, the actual
coupled equation for $\beta$ is exactly scalar on $r\le2r_c$
while the solution is smooth in this interval. These constants are independent of
$\lambda_0$.

Concavity of $\arctan$ gives $H_b(r)-\pi\ge\delta r$ for $r\le1$.
Since $\pi-Q_{\lambda_0}(r)\le2\lambda_0/r$, the condition
$\lambda_0\le\delta r_c^2/32$ gives
\[
 \beta_0\ge\pi+3\delta r_c/8
 \quad\hbox{on }[r_c/2,3r_c/2].
\]
Set $\eta_c=\delta r_c/8$, $B_c=16(\pi+1)$, $C_c=4B_c+2$.
On this annulus the explicit function
\[
 w(r,t)=\pi+2\eta_c-B_c(r/r_c-1)^2-C_c t/r_c^2
\]
is below the initial map, is negative at both endpoints, and satisfies
$w_t\le\mathcal T w$. Indeed its radial Laplacian is bounded below
by $-4B_c/r_c^2$ and the sine term by $-2/r_c^2$.
Comparison yields $\beta(r_c,t)\ge\pi+\eta_c$ for
$t\le\tau_\beta=\eta_c r_c^2/C_c$.

Now put
\[
 d_c=\frac{\delta}{16\sqrt{r_c}},\qquad
 \lambda_c(t)=(\sqrt{\lambda_0}-d_c t/8)^2,
 \qquad T_*=8\sqrt{\lambda_0}/d_c,
\]
and define $v=Q_{\lambda_c(t)}+d_c r^{3/2}$ on $[0,r_c]$.
The function $\sin z\cos z$ is $1$-Lipschitz, so
\[
 \mathcal T v\ge\frac{5d_c}{4}r^{-1/2},\qquad
 v_t=\frac{d_c}{2}\frac{r\sqrt{\lambda_c}}{\lambda_c^2+r^2}
                         \le\frac{d_c}{2}r^{-1/2}.
\]
Moreover $v(\cdot,0)\le\beta_0$ and
$v(r_c,t)\le\pi+\eta_c/2<\beta(r_c,t)$.
Choose $\lambda_0$ smaller so that
$T_*<\min(\tau_\alpha,\tau_\beta)$. This choice is compatible with the small-energy requirements below. Scalar comparison gives
\[
 \beta_r(0,t)\ge2/\lambda_c(t)\longrightarrow\infty
 \qquad(t\uparrow T_*).
\]
If the coupled solution loses smoothness earlier, finite-time
singularity is already established; otherwise this divergence
prevents continuation through $T_*$. The barrier is $C^1$ at zero
and smooth for $r>0$. Comparison remains valid: on every closed
time interval before $T_*$, the linearized zeroth-order coefficient
is nonpositive sufficiently near zero and bounded above elsewhere,
so a first negative minimum is excluded. This proves the first assertion of Theorem~\ref{thm:example}.

\subsection{The energy confines the two transitions}\label{sec:trap}

For the rest of the proof we consider the critical flow.
Its energy keeps each component close to a scalar sphere and keeps
the two transitions apart. The following argument proves this directly
from the initial data.

\paragraph{Compactness at the energy of one sphere.} If
$v:\R\to[0,5\pi/4]$ has endpoints $\pi,0$, then
\[
 \mathcal E(v)-4\pi
 =\pi\int_\R(v_s+\sin v)^2\,ds,
 \qquad \mathcal E(v)=\pi\int_\R(v_s^2+\sin^2v)\,ds.
\]
As $\mathcal E(v)\downarrow4\pi$, translation by any $\pi/2$ crossing
gives $\|v(\,\cdot+s_v)-Q\|_{H^1(\R)}\to0$, where $s_v$ is the crossing. Indeed, the first-order equation with
an $L^2$ error and Gronwall give convergence on fixed intervals. The
core energy tends to $4\pi\tanh K$ on $[-K,K]$. The remaining energy
is small. An excursion from an endpoint to $\pi/2$ costs at least
$2\pi$, by
\[
 \pi\int_a^b(v_s^2+\sin^2v)\,ds
 \ge2\pi|\cos v(b)-\cos v(a)|.
\]
Thus on each far tail $\sin^2v$ controls the squared distance to the
appropriate pole. This proves strong $H^1$ convergence, without
assuming a single crossing.

Write $v_\alpha,v_\beta$ for the following extensions to the line.
Extend $\alpha$ by $\pi$ for $s<0$, and $\beta$ by $\pi+\delta e^s$.
The latter extension costs at most $\pi\delta^2$. Since
$\mathfrak f\ge1$, each scalar energy is at least $4\pi$ (the
$\beta$ energy is in fact at least $4\pi+E_b$). It follows that
sufficiently small $E-8\pi$ and $\delta$ put both extensions in a
fixed arbitrarily small $H^1$ neighborhood of translates of $Q$.

Fix a nonnegative smooth $\xi\not\equiv0$ supported in $(-1,1)$ and
set $Z=\xi Q'$. Define $S_i$ by
\begin{equation}
 \int_\R[v_i(S_i+x)-Q(x)]Z(x)\,dx=0,\qquad i=\alpha,\beta.
 \label{eq:phase}
\end{equation}
The derivative is close to $\int\xi(Q')^2>0$, so these scales are
unique in the small neighborhood and depend smoothly on time.
Uniform closeness excludes other zeros away from that neighborhood.

Given a small tube size $\eta$ and large fixed $K_0$, sufficiently
small initial excess ensures, throughout the smooth lifetime,
\begin{equation}
 \|v_i(S_i+\cdot)-Q\|_{H^1(\R)}\le\eta,
 \qquad S_\alpha\ge K_0,\quad L=S_\beta-S_\alpha\ge K_0.
 \label{eq:tube}
\end{equation}
Indeed, the trace of the first extension at $s=0$ is $\pi$, so its
uniform approximation by $Q(-S_\alpha)$ forces $S_\alpha$ large.
Also $\alpha(S_\alpha)\ge\pi/2-C\eta$ and $\beta\ge\alpha$; the
uniform approximation $\beta(S_\alpha)=Q(-L)+O(\eta)$ excludes
$L\to-\infty$. For the separation
bound, argue by contradiction as the excess tends to zero: if the
relative shift did not tend to $+\infty$, then
\[
 \int(\mathfrak f(\alpha)-1)(\beta_s^2+\sin^2\beta)\,ds
\]
would have a positive limit on a fixed interval where the outer $Q$
lies in the nonconstant warping region. Ordering excludes reversal
of the two transitions; a bounded relative shift is impossible by
this positive interaction. The preceding scalar compactness proof
makes this a uniform statement for sufficiently small excess.

For clarity, write $z_b(r)=H_b(r)-\pi$. On $r\le1/2$ the second
angle is $Q_{\lambda_0}+z_b$ and $|z_b|+|(z_b)_s|\le C\delta r$.
The cross terms in the derivative and sine energies are bounded by
$C\delta\int_0^\infty e^{-s}\sech(s+\log\lambda_0)\,ds
=O(\delta\lambda_0(1+|\log\lambda_0|))$. For the sine term use
$|\sin^2(v+w)-\sin^2v-\sin^2w|\le4|\sin v\,\sin w|$.
The extra warping on the inner sphere costs
$O((\lambda_0/R_0)^2)$, since it vanishes for $r\le R_0$.
The warping factor multiplying the root energy equals $f_\pi$ outside $r\le C R_0$,
and the root energy in this latter disk is $O(\delta^2 R_0^2)$.
On $r\ge1/2$ the cutoff corrections to the two angles and their
$s$ derivatives are $O(R_0)$ and $O(\lambda_0)$, respectively.
Together these estimates give
\[
 E(u_0)=8\pi+f_\pi E_b+o(1)
 \quad(R_0\to0,\ \lambda_0/R_0\to0).
\]
Here the unweighted root energy is
$2\pi\int_\pi^{\pi+\delta}|\sin v|\,dv=E_b$. Thus $\delta$, then
$R_0$, then $\lambda_0$ can be chosen to give the required small
excess. Energy monotonicity yields \eqref{eq:tube} for as long as
the solution is smooth.

\subsection{The source near the two transitions}
\label{sec:source-bands}

We now analyze the critical flow. The parameter set selected by its
source must be determined from the current map, not prescribed in advance.
We first verify the normalization used in the construction.
At $p=(1,0)$, $\ell=1/8$, the
density on $B_{2\ell}(p)$ lies between $1024/1681$ and $1024/625$.
Integration over $B_\ell(p)$ gives $M\ge16\pi/1681$; using
$\int_{B_\ell(p)}|y-b|^2dy\ge\pi\ell^4/2$ and
$M\le(1024/625)4\pi\ell^2$ gives
$\sigma^2/\ell^2\ge625/(8\cdot1681)$. Weight one therefore holds on
an open parameter set whose physical supports avoid zero. This proves
$A_Q^{\rm raw}>0$ and will detect the outer sphere even when the inner
one has collapsed to an atom. To see finiteness, bounded density
excludes small radii. A packet of positive weight has a fixed
positive mass, so $|p|\le2\ell+C$ by the integrable tail of $\mu_Q$.
For large $\ell$ comparison of the variance with the moment about
zero bounds $\sigma^2/\ell^2$ by $C\log(2+\ell)/\ell^2$. This excludes
large radii and leaves a compact parameter box.

\paragraph{Localization at the two scales.}

Fix a time in the tube \eqref{eq:tube}, write
$R=e^{-S_\alpha}$, $\lambda=e^{-S_\beta}$, and put
$\varepsilon=\lambda/R=e^{-L}$. Below, a constant is uniform in
$S_\alpha,S_\beta,L$ unless a dependence is indicated. The fixed scan,
including its positive mass and variance onsets
$m=\varepsilon_{\rm pkt}/4$ and $\rho_-$, is
held fixed when the tube is tightened. Let $k_Q=2\pi/A_Q^{\rm raw}$.

\begin{lemma}[Localization in packet parameters]
\label{lem:source-boxes}
There are $0<c<C<\infty$ depending only on the fixed data such that,
for $\eta$ sufficiently small and $K_0$ sufficiently large in
\eqref{eq:tube}, every packet of positive weight belongs to one of
the disjoint boxes
\[
 \mathcal B_h=\{(p,\ell): |p|\le Ch,\ ch\le\ell\le Ch\},
 \qquad h\in\{\lambda,R\}.
\]
All corresponding cutoff supports lie in $D$, and $CR<r_*$.
For $h=\lambda,R$ define the parameter masses
\[
 \widehat A_h=k_Q\int_{\mathcal B_h}W_u(p,\ell)
                              \,\frac{dp\,d\ell}{\ell^3}.
\]
Then
\begin{equation}
 \widehat A_\lambda=2\pi+O(\eta+\varepsilon^2),\qquad
 \widehat A_R=k_Q A_{\rm dec}^{\rm raw}+O(\eta+\varepsilon),
 \qquad A_{\rm dec}^{\rm raw}>0,
 \label{eq:packet-masses}
\end{equation}
where $A_{\rm dec}^{\rm raw}$ is the unnormalized full-plane parameter
integral for the measure $\mu_Q+4\pi\delta_0$.
\end{lemma}

\begin{proof}
Let $d\mu_h=4h^2(h^2+|y|^2)^{-2}dy$. On $D$ one has
\begin{equation}
 \|\mu_u-(\mu_\lambda+\mu_R)|_D\|_{\rm TV}
       \le C(\eta+\varepsilon^2)=:\omega .
 \label{eq:packet-TV}
\end{equation}
Indeed the scalar energy in cylinder coordinates is
$\tfrac12(v_s^2+\sin^2v)\,ds\,d\theta$. Taking differences of squares
and applying Cauchy--Schwarz gives an $O(\eta)$ error for each scalar
measure. The warping factor changes by $O(\eta)$ in $L^\infty$ when
$\alpha$ is replaced by $Q(s-S_\alpha)$; the excess over one of this
last factor is supported in $|y|>R$, where $\mu_\lambda$ has mass
$4\pi\lambda^2/(R^2+\lambda^2)$. This proves \eqref{eq:packet-TV}.

We give the support argument, including packets not centered at the
origin. Scale $R$ to one, so the two reference radii are $1$ and
$\varepsilon$; write $p,\ell$ temporarily for the rescaled parameters.
Restricting the reference measures to the disk can only decrease their
contributions to the upper bounds below. Fix $\omega<m/8$. The reference
density is at most $C\varepsilon^{-2}$, so a packet with
$\ell<c\varepsilon$ cannot reach mass $m$. At scales $\ell\le c_1$
with $c_1$ small, the outer profile contributes at most $C\ell^2<m/4$.
Choose $K$ fixed so that the inner reference mass outside $B_{K\varepsilon}$
is less than $m/4$. Activity then forces
$|p|\le 2\ell+K\varepsilon$. In particular, if $\ell\ge K\varepsilon$,
the cutoff is supported in $B_{5\ell}$. The variance about the barycenter
is bounded by the second moment about zero, whether or not zero lies
in the packet. On this support the error in the latter moment is at most
$25\ell^2\omega$. Dividing by $M\ell^2\ge m\ell^2$ gives
\[
 \frac{\sigma^2}{\ell^2}\le C\left[
  \frac{\varepsilon^2}{\ell^2}\log(2+\ell/\varepsilon)
              +\ell^2+\omega\right].
\]
Here we used
$\int_{B_d}|y|^2d\mu_h\le Ch^2\log(2+d/h)$ for the inner profile
and the bounded density of the outer profile. Choosing a large
$C_1\ge K$, then a small $c_1$, and $\omega$ sufficiently small makes
the right side less than $\rho_-^2$ on
$C_1\varepsilon\le\ell\le c_1$.

At all scales, concentration of both reference masses in a fixed
$B_K$ up to a total tail less than $m/2$ gives $|p|\le2\ell+K$
for every active packet. For $\ell\ge K$ its support therefore lies
in $B_{5\ell}$, and the preceding second-moment argument gives
\[
 \frac{\sigma^2}{\ell^2}\le C\left[
 \frac{\log(2+\ell)+\varepsilon^2\log(2+\ell/\varepsilon)}
      {\ell^2}+\omega\right].
\]
This excludes $\ell\ge C_2$ uniformly for $0<\varepsilon\le1$:
$\varepsilon^2\log(2+\ell/\varepsilon)
\le\log(2+\ell)+\varepsilon^2\log(1/\varepsilon)$.
The two center bounds already obtained give the claimed boxes after
returning to physical units. Increasing $K_0$ ensures disjointness,
$3CR<1$, and $CR<r_*$. Thus neither the disk boundary nor the upper
cutoff removes a parameter in these boxes.

On a rescaled box $\mathcal B_h$ the weight is a uniformly Lipschitz
function of the cutoff moments of $1,Y,|Y|^2$, $Y=y/h$. In fact
the mass cutoff removes $M=0$ and the variance cutoff removes
the square-root singularity at zero; $\ell/h$ stays in a compact
subset of $(0,\infty)$. The parameter measure $dp\,d\ell/\ell^3$
is invariant under this rescaling, and each rescaled box has finite
parameter measure. The box constants can be enlarged to contain the
whole active scan of the reference $\mu_Q$ at scale one, and also of
$\mu_Q+4\pi\delta_0$: for the latter, at small $\ell$ a packet reaching
mass $m$ must see the atom, hence $|p|\le2\ell$, and its variance ratio
squared is $O(\ell^2)$ by comparison with zero. The large-radius
argument above is unchanged.

In the inner chart all moment tests have bounded support, on which
the extra outer reference mass is $O(\varepsilon^2)$. The difference
from the isolated-profile moments is thus $O(\omega+\varepsilon^2)$.
In the outer chart, each moment test $P$ has uniformly bounded
Lipschitz norm and support inside the rescaled disk. Since
$\int_{\R^2}|Y|\,d\mu_\varepsilon(Y)=O(\varepsilon)$,
\[
 \left|\int P\,d\mu_\varepsilon-4\pi P(0)\right|
 \le \operatorname{Lip}(P)\int |Y|\,d\mu_\varepsilon(Y)
 \le C\varepsilon .
\]
Applying moment stability and integrating over the boxes proves
\eqref{eq:packet-masses}. The strictly accepted off-origin parameter
block for $\mu_Q$ at the beginning of this subsection is disjoint from the atom on its
cutoff supports and therefore also proves $A_{\rm dec}^{\rm raw}>0$.
\end{proof}

\paragraph{The heat on the disk.}
Split the actual source as $\nu=\nu_\lambda^D+\nu_R^D$ by restricting
the defining integral to the two boxes. Write $\nu_h^E$ for the same
integral with the plane heat kernel, retaining the very same weights
and barycenters. All four functions are radial: rotation preserves
the packet weight of a radial energy measure, rotates its barycenter,
and preserves the corresponding heat kernel. This observation holds
at every smooth time, without any limit of the profiles. Domination
by the plane heat kernel and $|b_{p,\ell}|\le3Ch$ give
\begin{equation}
 0\le\nu_h^D(y)\le\nu_h^E(y)
       \le C h^{-2}e^{-c|y|^2/h^2}\quad(y\in D),
 \qquad \int_{\R^2}\nu_h^E=\widehat A_h .
 \label{eq:band-gaussians}
\end{equation}
The Gaussian upper bound holds on the entire plane for $\nu_h^E$.
Let $A_h^D=\int_D\nu_h^D$. Uniformly in either box,
\begin{equation}
 0\le \widehat A_h-A_h^D\le C e^{-c/h^2}.
 \label{eq:killed-mass}
\end{equation}
Indeed, the boundary values of the free kernel, at every preceding
heat time $0<\tau\le\ell^2$, are bounded by $Ch^{-2}e^{-c/h^2}$
since $|b_{p,\ell}|\le1/2$. The caloric difference from the killed kernel
has zero initial data and those boundary values. The maximum principle
gives this same bound inside $D$ (equivalently, first regularize the
unit point mass). Integrate over the finite rescaled parameter box
and over the disk, and add the Gaussian tail outside the disk; the
polynomial prefactor is absorbed by decreasing $c$. In particular
$G_D(\nu_h^E|_D-\nu_h^D)$ is bounded by $Ce^{-c/h^2}$, since
$G_D1=(1-|y|^2)/4$.

For definiteness use the packet residues
$A_{\rm in}=\widehat A_\lambda$, $A_{\rm out}=\widehat A_R$.
Equations \eqref{eq:packet-masses}--\eqref{eq:killed-mass} imply
$A_\lambda^D+A_R^D>2\pi+c_0$ for some $c_0>0$, after the tube choices;
no equality of the two limiting residues is used. The only disk mass
defect that can multiply $L$ in the middle is the \emph{inner} one:
as $\lambda=Re^{-L}\le e^{-L}$, for any fixed $c'>0$,
\begin{equation}
  (1+L)|A_{\rm in}-A_\lambda^D|\le C_{c'}e^{-c'L}.
 \label{eq:inner-heat-loss}
\end{equation}
In particular, using the true disk inner mass for the residue instead
would change none of the later conclusions.

\paragraph{The radial coefficient.}
Let $m(r)=\int_{B_r}\nu$. The Dirichlet condition and radial Poisson
equation give the exact identities
\[
 \phi(r)=\frac1{2\pi}\int_r^1 m(z)\,\frac{dz}{z},\qquad
 \partial_s\log q=2-\frac{m(e^{-s})}{\pi},\qquad
 \partial_s^2\log q=2r^2\nu(r)\ge0 .
\]
Put $q_\alpha=q(S_\alpha)$, $q_\beta=q(S_\beta)$,
$q_m=\inf_{s\ge0}q(s)$, and $d=A_{\rm in}/\pi-2$.
Directly for the disk source, \eqref{eq:band-gaussians} gives
\begin{equation}
 |m(r)-A_\lambda^D|
 \le C e^{-c(r/\lambda)^2}+C(r/R)^2
 \qquad(\lambda\le r\le R).
 \label{eq:disk-middle-mass}
\end{equation}
The difference on the left, before taking its absolute value, is the
outer mass inside $B_r$ minus the inner mass outside $B_r$. Both
error bounds are integrable against $dr/r$ with uniform integrals.
Since $s=-\log r$, integration with radius increasing gives
$\log(q(-\log r)/q_\beta)=
\int_\lambda^r(m(z)/\pi-2)\,dz/z$, with a \emph{plus} sign
for the mass term. Consequently on the closed interval
$\lambda\le r\le R$,
\[
 \log\frac{q(-\log r)}{q_\beta}
        =d\log(r/\lambda)+O(1)
\]
by \eqref{eq:inner-heat-loss}. For $r\le\lambda$, one has
$m(r)\le C(r/\lambda)^2$; hence
$0\le\phi(r)-\phi(\lambda)\le C$ and
$q(-\log r)\ge c q_\beta(\lambda/r)^2$.
For $R\le r\le1$ the mass missing from $A_\lambda^D+A_R^D$ is at most
$Ce^{-c(r/R)^2}$, again integrable against $dr/r$. Thus
\[
 \log\frac{q(-\log r)}{q_\alpha}
 =\left(\frac{A_\lambda^D+A_R^D}{\pi}-2\right)\log(r/R)+O(1),
\]
whose logarithmic coefficient is positive. On the middle interval
an affine function of $\log r$ takes its minimum at an endpoint.
Together with $q(0)=1$ in the $s$ coordinate, the three estimates prove
\begin{equation}
 \log\frac{q_\beta}{q_\alpha}=-dL+O(1),\qquad
 q_m\asymp\min(q_\alpha,q_\beta),\qquad q_m\le1 .
 \label{eq:radial-ratios}
\end{equation}
Also $|\partial_s\log q|\le C$ by the total mass bound, giving
two-sided comparison on fixed collars of either core. Given $d_0>0$,
choose $\eta$ small and $K_0>2$ large enough that $|d|\le d_0$ by
\eqref{eq:packet-masses}. Without assigning a sign to $d$ we obtain
\begin{equation}
 \frac{q(S_\beta+x)}{q_\beta}\ge
 \begin{cases}
 c e^{-d_0|x|},&-L/2-1\le x\le0,\\
 c,&x\ge0.
 \end{cases}
 \label{eq:inner-q-weight}
\end{equation}

We record the version used for metric convergence. For $K\ge1$,
$L\ge2\log K$ and $r_m=\sqrt{\lambda R}$, integration of
\eqref{eq:disk-middle-mass} from $r_m$ to any point of the trimmed
interval gives
\begin{equation}
 \sup_{K\lambda\le r\le R/K}
 \left|\log\frac{q(-\log r)}{q(-\log r_m)}\right|
 \le |d|L+CK^{-2}+Ce^{-cK^2}+CL e^{-c/\lambda^2}.
 \label{eq:trimmed-packet-estimate}
\end{equation}
Indeed the outer error integrates to at most
$C\int_0^{1/K}z\,dz$, and the inner error to at most
$C\int_K^\infty e^{-cz^2}\,dz/z$. The last term is controlled by
\eqref{eq:inner-heat-loss}; no estimate here multiplies an outer heat
loss or a fixed profile error by the neck length in place of $d$.

\paragraph{The error in the inner source mass.}
\label{sec:residuetime}

The packet estimate by itself gives a small number $|d|$, but not
a bound on $L|d|$. The following estimate replaces that small number
by a quantity measured on the same time slice.

\begin{lemma}[Residue--tension estimate]\label{lem:residue}
At every smooth time in the tube \eqref{eq:tube},
\begin{equation}
 |d|\le C\sqrt{\Diss/q_\beta}+Ce^{-cL}.
 \label{eq:residue-tension}
\end{equation}
The constants do not depend on $S_\alpha$ or $S_\beta$.
\end{lemma}

\begin{proof}
Write $x=s-S_\beta$ and $h(x)=\beta(S_\beta+x)-Q(x)$; below, the
untranslated symbols for angles and tensions in this chart refer to
their values at $S_\beta+x$. On $x\ge -L/2-1$ the outer angle obeys
$0\le\alpha\le C\eta+Ce^{-L/2}$ by \eqref{eq:tube} and the
one-dimensional Sobolev inequality. Here $\mathfrak f=1$ and both
equations have the scalar tension $F_v=v''-\sin v\cos v$.

\emph{An inverse with a one-sided weight.}
Put $H=-\partial_x^2+1-2\sech^2x$ and $B=\partial_x+\tanh x$.
Then $H=B^*B$, $BB^*=-\partial_x^2+1$, and $\ker B=\R Q'$,
with $Q'=-\sech x$. In particular $H$ has gap at least $1$ on
$(Q')^\perp$: the range of $B^*$ equals $(\ker B)^\perp$
because $BB^*\ge1$, and the nonzero spectra of $B^*B$ and $BB^*$
agree. Controlling the kernel coefficient by the phase test $Z$ in
\eqref{eq:phase}, for which $\langle Q',Z\rangle>0$, and then using
$v''=(1-2\sech^2x)v-Hv$ gives
\begin{equation}
 \|v\|_{H^2(\R)}\le C\bigl(\|Hv\|_2+|\langle v,Z\rangle|\bigr).
 \label{eq:unweighted-inverse}
\end{equation}
Choose smooth $\psi$ equal to $x$ for $x\le-1$ and to zero for
$x\ge0$, and let $w=e^{\gamma\psi}$ with small fixed $\gamma>0$.
Conjugation gives
\[
 wHw^{-1}=H+2\gamma\psi'\partial_x+
              \gamma\psi''-\gamma^2(\psi')^2 .
\]
The added terms are bounded by $C\gamma$ from $H^2$ to $L^2$.
Also
$|\langle wv,Z\rangle-\langle v,Z\rangle|
\le C\gamma\|wv\|_2$, since $Z$ is compactly supported. Apply
\eqref{eq:unweighted-inverse} to $wv$ and absorb:
\begin{equation}
 \|wv\|_{H^2}\le C\bigl(\|wHv\|_2+|\langle v,Z\rangle|\bigr).
 \label{eq:weighted-inverse}
\end{equation}
For $-\partial_x^2+P$ with $0<p_0\le P\le1$ the same estimate
holds without the phase term, by the coercive quadratic form and
the same second-derivative estimate. Choose $\gamma$ for these two
inverses first, and shrink the tube so that $d_0<2\gamma$ in
\eqref{eq:inner-q-weight}.

\emph{Localization to the entire inner disk.}
Take $0\le\chi_L\le1$ vanishing on $x\le-L/2-1$ and equal to one on
$x\ge-L/2$, with uniformly bounded derivatives. Extend the cutoff
functions by zero to the whole line. The cutoff is one on
$\operatorname{supp} Z$, so $\langle\chi_Lh,Z\rangle=0$.
Since $Q''=\sin Q\cos Q$ and $\cos(2Q)=1-2\sech^2x$, Taylor's
formula gives
\[
 H(\chi_Lh)=-\chi_LF_\beta+
 O(\eta\chi_L|h|)-2\chi_L'h'-\chi_L''h.
\]
After multiplication by $w$, the last two terms have $L^2$ norm
$O(e^{-\gamma L/2})$ by the uniform $H^1$ bound. Moreover
$w^2\le Cq(S_\beta+x)/q_\beta$ throughout the cutoff support,
by \eqref{eq:inner-q-weight}. Thus for $i=\alpha,\beta$,
\[
 \|w\chi_LF_i\|_2\le C\sqrt{\Diss/q_\beta}.
\]
Apply \eqref{eq:weighted-inverse} and absorb $C\eta\|w\chi_Lh\|_2$.
For $\chi_L\alpha$ the identical commutator calculation uses
$P=\sin\alpha\cos\alpha/\alpha$, defined as $1$ at zero and
extended off the cutoff region with a fixed positive lower bound.
Here the transition-strip $H^1$ norm of $\alpha$ is uniformly
bounded by the tube estimate and the outer $Q$ tail. Consequently,
for each fixed $C_1$ containing all the dilated inner packet
supports in $B_{C_1}$, taking $L/2>\log C_1+1$ gives
\[
 \|h\|_{H^1(x\ge-\log C_1)}
   +\|\alpha\|_{H^1(x\ge-\log C_1)}
 \le C\bigl(\sqrt{\Diss/q_\beta}+e^{-cL}\bigr).
\]
The weight is one towards $+\infty$: this controls the disk all
the way to its center, not merely an annulus around $r=\lambda$.

On this disk the two scalar energy densities add with unit weights.
By Cauchy--Schwarz for their derivative and sine squares, the
total variation distance of the rescaled measure from $\mu_Q$ is
bounded by $C\|h\|_{H^1}+C\|\alpha\|_{H^1}^2$
(all norms on the preceding half-line). The normalized mass,
first moment and second moment therefore differ by at most the same
bound times a fixed constant on the inner parameter box. Choose the
box to contain the active parameters of both the actual inner scan
and the isolated $Q$ scan. The smooth combined weight is Lipschitz
there, the normalized parameter volume is bounded, and the isolated
mass is $2\pi$. Thus, with $z=\Diss/q_\beta$,
$|d|\le C(\sqrt z+z+e^{-cL})$. For $z\le1$ absorb $z$ into
$\sqrt z$; for $z>1$ use the uniform packet bound on $d$.
\end{proof}

\begin{proposition}[Scale motion at critical calibration]
\label{prop:absorption}
On the smooth trajectory,
\begin{equation}
 q_\beta\le Cq_m+CL^2\Diss,\qquad
 |\dot S_\beta|\le C\sqrt{q_m\Diss}+CL\Diss .
 \label{eq:refined-motion}
\end{equation}
\end{proposition}
\begin{proof}
By \eqref{eq:radial-ratios}, \eqref{eq:residue-tension} and boundedness
of $Le^{-cL}$, there are fixed $C_0,C_2$ such that
\[
 \log(q_\beta/q_\alpha)\le C_0+C_2 L\sqrt{\Diss/q_\beta}.
\]
If the left side is at most $C_0+1$, the minimum comparison in
\eqref{eq:radial-ratios} gives $q_\beta\le Cq_m$. Otherwise
$1\le C_2L\sqrt{\Diss/q_\beta}$, so $q_\beta\le C_2^2L^2\Diss$.
To estimate motion, differentiate \eqref{eq:phase}. The denominator
$\int\beta_s(S_\beta+x)Z(x)\,dx$ is bounded below, by the tube
and $\langle Q',Z\rangle>0$. On $\operatorname{supp}Z$,
$q(S_\beta+x)\le Cq_\beta$ (the logarithmic derivative of $q$
is uniformly bounded by the packet mass bound). Thus
\[
 |\dot S_\beta|^2
 \le C\left(\int\frac{\beta_t^2}{q}\,ds\right)
       \left(\int q(S_\beta+x) Z(x)^2\,dx\right)
 \le C q_\beta\Diss .
\]
Here $\beta_t=qF_\beta$ and $\mathfrak f\ge1$. Substituting the
first estimate gives the second. No sign of $d$ is needed.
\end{proof}

\subsection{Global existence and simultaneous small tension}
\label{sec:continuationclock}

\begin{proposition}[No finite-time singularity]\label{prop:continuation}
The smooth critical flow exists at all finite times.
\end{proposition}
\begin{proof}
\emph{The two positions stay bounded on finite intervals.}
Write $Y=1+S_\beta$, so $1\le Y$ and $L\le Y$. By
\eqref{eq:refined-motion} and $q_m\le1$,
\begin{equation}
 |\partial_t\log Y|\le C\frac{\sqrt{q_m\Diss}}Y+C\Diss
 \le C\sqrt{\Diss}+C\Diss .
 \label{eq:logY}
\end{equation}
The right side is integrable up to any finite maximal time $T$,
since $\int_0^T\Diss\le E(u_0)$ and Cauchy--Schwarz applies on
bounded intervals. Therefore $S_\alpha\le S_\beta$ remain bounded
above, and $R,\lambda$ have positive lower bounds. All active
heat times are now bounded below. Differentiating the disk heat
kernel in its observation variable and solving the fixed Dirichlet
Poisson problem gives
\[
 c_T\le a\le1,\qquad
 \sup_{0\le t<T}\big(\|\phi(t)\|_{C^k(\overline D)}
                 +\|a(t)\|_{C^k(\overline D)}\big)\le C_{k,T}
\]
for every $k$. No time regularity of $a$ or derivative bound for
the unknown map is used in this step.

\emph{A regular equation on a cap containing the pole.}
The tube and boundedness of $S_\beta$ give a fixed $r_c>0$ such
that both angles lie in $[0,\pi/4]$ on $r\le r_c$ for $t<T$.
The warping is constant there. For each angle $v$ use the stationary
comparison solution $P_\rho(r)=2\arctan(r/\rho)$ of the scalar
equation. One $\rho>0$ suffices for both angles: require
$\rho\sup_{0<r\le r_c}\tan(v(r,0)/2)/r\le1$ and
$P_\rho(r_c)\ge\pi/4$. The supremum is finite by initial
smoothness. Thus $0\le v\le P_\rho\le 2r/\rho$ throughout the
cap. For clarity, comparison can be made without a singular
maximum principle: on every compact smooth time interval the lifts
$w=v/r$ and $P_\rho/r$ are smooth radial functions of $y\in\R^4$,
$r=|y|$, and satisfy the same equation
\begin{equation}
 w_t=a(\Delta_{\R^4}w+N(r,w)),\qquad
 N(r,w)=\frac{rw-\sin(rw)\cos(rw)}{r^3}.
 \label{eq:polelift}
\end{equation}
Smooth corotationality at the north pole means that $v/r$ is
smooth and even in $r$ (use the tangent-plane coordinates of each
sphere near the pole). Similarly the radial coefficient $a$ has a
smooth even lift, with uniform spatial bounds inherited from the
disk. The quotient $N$ is smooth in $(r^2,w)$ on bounded sets
by its Taylor expansion. Scalar comparison for \eqref{eq:polelift}
gives the claimed barrier and hence a uniform bound on $w$.

For completeness a local derivative bound does not require $a_t$.
Differentiating \eqref{eq:polelift} spatially, the term
$2(Da\cdot Dw)\Delta w$ is absorbed by part of
$-2a|D^2w|^2$. Hence
\[
 (\partial_t-a\Delta)|Dw|^2
 \le -a|D^2w|^2+C(1+|Dw|^2),\qquad
 (\partial_t-a\Delta)w^2\le C-2a|Dw|^2 .
\]
For a cutoff $\zeta$ compactly supported in $B_{r_c}$ and equal to
one on $B_{r_c/2}$, apply the maximum principle to
$\zeta^2|Dw|^2+B w^2$. The cutoff cross term is bounded by
$\tfrac12a\zeta^2|D^2w|^2+C|Dw|^2$, and fixed large $B$ absorbs
the latter. The initial and lateral values are bounded. This gives
uniform bounds for both $w,Dw$ on the smaller cap.

\emph{A boundary estimate on the remaining annulus.}
Let $U=(\alpha,\beta)$ be the curve in the smooth meridional
metric $d\alpha^2+\mathfrak f(\alpha)d\beta^2$, understood as a
metric on the lifted angle plane (extend the warping constantly
near its endpoints). Put
\[
 V(U)=\tfrac12(\sin^2\alpha+\mathfrak f(\alpha)\sin^2\beta),
 \qquad X=U_s,\qquad A=\nabla_s X .
\]
Then $U_t=q(A-\nabla V)$. On the fixed interval
$0\le s\le s_c=-\log(r_c/3)$, $q,q^{-1},q_s$ are uniformly
bounded. There is an exact scalar inequality up to the disk
boundary. Namely for $J=|X|^2/2-V(U)$,
\begin{equation}
 (\partial_t-q\partial_s^2-q_s\partial_s)J
     =q\bigl(|\nabla V|^2-|A|^2\bigr)\le C .
 \label{eq:annular-gradient}
\end{equation}
Indeed $\nabla_tX=\nabla_s U_t$; subtraction of the equation
for $V(U)$ cancels both Hessian terms. At $s=0$ the fixed
Dirichlet data give $U_t=0$, hence $A=\nabla V$ and
$J_s=\langle A-\nabla V,X\rangle=0$. At $s=s_c$ the cap
estimate bounds $J$. The scalar maximum principle with this
homogeneous Neumann condition at one end bounds $J$ from above.
Thus $|U_s|$ is uniformly bounded on the entire fixed annulus.

\emph{Closing the continuation argument.}
In the cap use \eqref{eq:polelift}; on the fixed annulus use the
angle equations in coordinates. Their inhomogeneous terms are
now bounded. The interior parabolic $W^{2,1}_p$ estimate for
uniformly elliptic, spatially Lipschitz leading coefficients,
measurable in time (a special case of \cite{KimKrylov}), gives
uniform local $W^{2,1}_p$ norms on strips away from the initial
time, for $p>6$. These estimates hold up to $T$ from below, by
the forward estimates on truncated time intervals. At $s=0$,
subtract the constant boundary values, extend the components and
the bounded right sides oddly, and extend $q$ evenly; the leading
coefficient remains spatially Lipschitz. The cap and annulus
interiors overlap on $r_c/3<r<r_c/2$, giving global bounds.
Parabolic Sobolev embedding
therefore gives the lifted functions and their first spatial
derivatives a uniform spatial and temporal H\"older modulus.
This also holds for the original map in disk coordinates:
in a cap chart its horizontal tangent-plane coordinates are
$y\,\sin(v)/r$ with $y\in\R^2$, and its vertical coordinate is
$\cos(v)$, both smooth functions of $y$ and the restricted lift
$w(y,0,0)$. The positive packet-radius cutoff makes $u\mapsto a_u$
locally Lipschitz from $C^1(\overline D)$ to each spatial $C^k$;
this follows by differentiating the cutoff moments, positive-time
heat kernel and Poisson operator. It gives a uniform temporal
H\"older modulus for $a$ as well. Interior and fixed-Dirichlet
boundary Schauder estimates now give $C^{2+\theta,1+\theta/2}$
control away from the initial time for some $\theta>0$.
In bootstrapping time derivatives there is no loss of a spatial
derivative in a packet moment: for each fixed smooth moment test $h_0$,
\[
 \frac{d}{dt}\int_D h_0\,d\mu_u
 =-\int_D\langle h_0\tau(u)+du(\nabla h_0),u_t\rangle\,dx
\]
by fixed Dirichlet data. This gives a uniform H\"older bound on
$a_t$ from the current bounds on $D^2u,u_t$; derivatives of the
tests are uniform on the truncated parameters. Differentiate the
equation in time (the differentiated boundary data vanish), use
Schauder estimates on a smaller time strip, and recover the higher
spatial derivatives from the elliptic equation for each time slice.
Iteration, with the displayed moment identity differentiated at
each stage, gives smooth terminal data at $T$ with the Dirichlet
compatibility conditions.

We recall why the local theory used to restart is valid for this
nonlocal coefficient. On a $C^1$ ball about any smooth datum,
the uniform slope bound makes the packet mass $O(\ell^2)$, so
the positive onset removes all sufficiently small $\ell$ uniformly.
On the remaining compact parameters the cutoff formula is smooth
in the moments, hence $a[u]$ is positive and locally Lipschitz in
$C^1$ with values in spatial $C^k$. Embed the target isometrically
and use the usual smooth tubular extension of the extrinsic map
equation; its lower-order term is smooth in $(u,Du)$ and the
target is preserved by the squared-distance maximum principle.
Freeze the leading coefficient at the initial datum. The linear
Dirichlet heat estimate bounds the solution with zero initial and
boundary values in $W^{2,1}_p$ by the $L^p$ norm of the forcing.
On bounded $W^{2,1}_p$ balls with fixed initial and boundary data,
the zero-initial-trace embedding bounds differences in
$L^\infty_t C^1_x$ by $C T_0^\gamma$ times their $W^{2,1}_p$
distance, for short time $T_0$ and some $\gamma>0$ ($p>6$).
In particular $(a[u]-a[u(0)])\Delta u$ and the lower-order term
give a contraction for small $T_0$: the factor multiplying
$\Delta u$ is uniformly small and its difference is
$O(T_0^\gamma)$ times the distance. One may center the ball at
the constant-in-time smooth datum, whose residual has $L^p$ norm
$O(T_0^{1/p})$. This proves local existence and uniqueness;
Schauder regularity recovers the smooth solution for compatible
data, and equivariance follows from uniqueness. Applied at $T$,
this contradicts finite maximality.
\end{proof}

\begin{proposition}[Simultaneous small tension at the two scales]\label{prop:clock}
For this global trajectory,
\begin{equation}
 \int_0^\infty\frac{q_m(t)}{(1+S_\beta(t))^2}\,dt=\infty .
 \label{eq:strong-clock}
\end{equation}
There is one sequence $t_j\to\infty$ such that
\begin{equation}
 (1+S_\beta(t_j))^2\,\frac{\Diss(t_j)}{q_m(t_j)}\longrightarrow0 .
 \label{eq:strong-slices}
\end{equation}
Both tensions tend to zero on the whole cylinder in $L^2(ds)$
on this same sequence.
\end{proposition}
\begin{proof}
Suppose the integral in \eqref{eq:strong-clock} were finite. The
first inequality in \eqref{eq:logY} gives
\[
 \int_0^\infty|\partial_t\log Y|\,dt
 \le C\left(\int_0^\infty\frac{q_m}{Y^2}\,dt\right)^{1/2}
       \left(\int_0^\infty\Diss\,dt\right)^{1/2}
        +C\int_0^\infty\Diss\,dt<\infty .
\]
Then $Y$ is bounded. Both packet radii have uniform positive lower
bounds, so the heat and Poisson estimates above give
$a\ge c>0$ for all time. Since $q=a/r^2\ge a$ on $D$,
$q_m/Y^2$ is bounded away from zero, a contradiction.
For the infinite measure $d\omega=(q_m/Y^2)\,dt$,
the nonnegative function $Y^2\Diss/q_m$ has integral
$\int_0^\infty\Diss\,dt<\infty$; every finite time interval has
finite $\omega$-measure since $q_m/Y^2\le1$. Thus for every $j$
there are arbitrarily large times where $Y^2\Diss/q_m<1/j$.
Choosing one increasing sequence proves \eqref{eq:strong-slices}.
Finally
\[
 \int_0^\infty(F_\alpha^2+F_\beta^2)(s,t_j)\,ds
 \le\frac{\Diss(t_j)}{2\pi q_m(t_j)}\longrightarrow0 .
\]
Neither separation in the limit nor convergence to a stationary
profile was assumed in selecting these times.
\end{proof}

\subsection{Two distinct harmonic bubbles form at infinite time}
\label{sec:limits}

Write $v_j=(\alpha,\beta)(\cdot,t_j)$ at the times given by
\eqref{eq:strong-slices}. On each bounded interval of the fixed or a
translated cylinder, the energy bounds $v_j'$ in $L^2$. The equations
for $F_\alpha,F_\beta$ then bound $v_j''$ first in $L^1$ (the
quadratic terms are in $L^1$), hence $v_j'$ in $L^\infty$, and
finally $v_j''$ in $L^2$. Compactness in $C^1$ and the vanishing
$L^2$ tensions give a stationary limit after passage to a subsequence.
This applies up to $s=0$ in the disk chart. All limits have finite
energy by lower semicontinuity.

We use the tube \eqref{eq:tube} also in the uniform norm, by the
one-dimensional Sobolev inequality. In particular, on either side
of a translated transition, outside a fixed large core, the angle
is confined to a small neighborhood of the corresponding pole. In a
stationary limit on an infinite tail the angle actually tends to that
pole: its signed distance from the pole has finite $H^1$ norm there,
since its square is comparable to $\sin^2 v$. This observation
identifies the endpoints below without assuming monotonicity.

\paragraph{The two transitions cannot belong to one harmonic sphere.}

Suppose a finite-energy stationary pair on $\R$ carries both
transitions, with endpoints $(\pi,\pi)$ and $(0,0)$ in order of
increasing $s$. Set
\[
 w=\alpha_s+\sin\alpha,\qquad
 P=\tfrac12\mathfrak f'(\alpha)(\beta_s^2+\sin^2\beta)\ge0.
\]
Then $w_s-\cos\alpha\,w=P$. The integrating factor
$I(s)=\exp(-\int_0^s\cos\alpha(\sigma)\,d\sigma)$ decays exponentially
at both ends. Since $w\in L^2$, one may integrate $(Iw)'=IP$ along
endpoints tending to $\pm\infty$ with $w$ tending to zero there.
It follows that $P\equiv0$. Wherever $\alpha$ passes through
$(\pi/2,3\pi/4)$ this forces $\beta_s=\sin\beta=0$.
Uniqueness for the stationary $\beta$ equation with this fixed
$\alpha$ forces $\beta$ to be constant, contradicting its endpoints.

\paragraph{The outer transition cannot remain in the disk.}

Suppose the outer center stays bounded on a subsequence. The resulting
stationary disk limit has $\alpha(0)=\pi$, $\alpha(+\infty)=0$
and $0\le\alpha\le\pi$. It cannot touch either pole at an interior
point: at such an extremum its derivative would vanish, and uniqueness
in the flat warping region would make it constant. On its final tail
$\mathfrak f=1$; finite energy gives the zero first integral and
$\alpha_s=-\sin\alpha$. More generally
\[
 J_\alpha:=\alpha_s^2-\sin^2\alpha,\qquad
 (J_\alpha)_s=\mathfrak f'(\alpha)\alpha_s
                              (\beta_s^2+\sin^2\beta).
\]
Continue backwards on the maximal interval where $\alpha_s<0$.
There $J_\alpha\ge0$, since $J_\alpha=0$ on the final tail.
An interior first zero of $\alpha_s$ would yield
$J_\alpha=-\sin^2\alpha<0$. Thus $\alpha_s<0$ throughout,
and $h:=-\alpha_s(0)>0$ by uniqueness at the outer pole. Similarly
\[
 J_\beta:=\beta_s^2-\sin^2\beta,\qquad
 (J_\beta)_s=-2\frac{\mathfrak f'(\alpha)}{\mathfrak f(\alpha)}
                         \alpha_s\beta_s^2\ge0.
\]
Since $J_\beta\in L^1$ and is nondecreasing, $J_\beta\le0$.
If the inner transition is also retained, $\beta$ must pass from
$b>\pi$ to $0$. At a crossing of $\pi$ one would have $\beta_s=0$;
uniqueness would then give $\beta\equiv\pi$, a contradiction.

In the remaining case the inner center escapes. The tube confines
$z:=\beta-\pi$ to a small neighborhood of zero on the whole limiting
half-line, with $z(0)=\delta$ and $z(+\infty)=0$. It obeys
\[
 (\mathfrak f z_s)_s=\mathfrak f\sin z\cos z .
\]
A negative minimum is impossible and an interior zero would force
$z\equiv0$. Thus $z>0$. The flux $\mathfrak f z_s$ is strictly
increasing and must stay negative, since $z$ tends to zero. In
particular $0<z\le\delta$ and, by $J_\beta\le0$,
$|z_s|\le\sin z\le z$.

Let $s_+<s_-$ be the crossings of $3\pi/4$ and $\pi/2$ by $\alpha$.
Since $-\alpha_s\ge\sin\alpha$, their distance is bounded by
$\int_{\pi/2}^{3\pi/4}(\sin u)^{-1}\,du$. For $s\ge s_-$ we have
$z_s=-\sin z$. Increasing flux and $|z_s|\le z$ therefore give,
with $v=z(s_+)$,
\[
 -z_s(s_+)\ge f_\pi^{-1}\sin z(s_-)\ge c v .
\]
On $[0,s_+]$ the warping is constant. The first integral for $\alpha$
gives the one-sided estimate
\[
 s_+=\int_0^{\pi/4}\frac{du}{\sqrt{\sin^2u+h^2}}
 \ge\operatorname{arsinh}\!\left(\frac{\pi}{4h}\right),
 \qquad e^{-s_+}\le C h .
\]
For $z$ put $k^2=\sin^2v-z_s(s_+)^2\ge0$. Its first integral gives
\[
 s_+=\int_v^\delta\frac{du}{\sqrt{\sin^2u-k^2}}
       \le\log(\delta/v)+C .
\]
Indeed the denominator is bounded below by $c' u$ for
$v\le u\le\delta$, using $-z_s(s_+)\ge cv$ and the monotonicity of
$\sin u$ in this range. Subtraction of $1/u$ from the integrand
leaves at most $C(u+k^2/u^3)$, whose integral is bounded uniformly
since $k^2\le v^2$. Thus $v\le C\delta e^{-s_+}$. Integrating the
$J_\alpha$ identity across the warping band now yields
\[
 h^2=\int_{s_+}^{s_-}\mathfrak f'(\alpha)(-\alpha_s)
                             (z_s^2+\sin^2z)\,ds
 \le 2(f_\pi-1)v^2\le C\epsilon_w\delta^2 h^2.
\]
The constants may be taken uniform for small $\epsilon_w,\delta$.
This contradicts $h>0$ and excludes the second case.

\paragraph{Identification of the root and the two spheres.}

The disk obstruction implies $S_\alpha(t_j)\to\infty$: otherwise a
subsequence retains the outer transition by \eqref{eq:tube}, while
the inner center is either bounded or escapes. If $L(t_j)$ had a
bounded subsequence, translation by $S_\alpha(t_j)$ would likewise
give the excluded merged stationary pair. Thus $L(t_j)\to\infty$.

The root stays near $(\pi,\pi)$ throughout its limiting half-line,
where the warping is constant. The zero scalar first integrals and
the boundary data give $\alpha\equiv\pi$ and
$\beta(s)=H_b(e^{-s})=\pi+2\arctan(c_b e^{-s})$. At the outer vertex
$z=\beta-\pi$ is small on all of $\R$ and lies in $H^1(\R)$.
Testing $(\mathfrak f z_s)_s=\mathfrak f\sin z\cos z$ with $z$,
using cutoffs tending to one, gives
$\int\mathfrak f(z_s^2+z\sin z\cos z)=0$, hence $\beta\equiv\pi$.
The first angle then solves the scalar stationary equation between
$\pi$ and $0$. Its zero first integral gives a translate of $Q$, and
\eqref{eq:phase} fixes that translate to be $Q$ itself. At the inner
vertex $\alpha$ is small on all of $\R$, so $\mathfrak f=1$ there;
the same test gives $\alpha\equiv0$, and the scalar first integral
and phase give $\beta=Q$. The limits are unique, so convergence in
each fixed chart holds on the selected sequence without further
extraction.

For fixed large $J$, the complementary intervals $I_J^0,I_J^1,I_J^2$
of Theorem~\ref{thm:example} lie in fixed small neighborhoods of
$(\pi,\pi),(0,\pi),(0,0)$, respectively, by \eqref{eq:tube}.
In each the warping is constant. For either component let $w$ be
its angle minus the indicated pole. There $F=w_{ss}-\sin w\cos w$
and $w\sin w\cos w\ge c w^2$. Testing by $w$ and using Young's
inequality gives
\begin{equation}
 \int_I(w_s^2+w^2)\,ds
 \le C\int_I F^2\,ds+C\big|[ww_s]_{\partial I}\big|.
 \label{eq:edge-coercivity}
\end{equation}
On the terminal half-line truncate along endpoints with $ww_s\to0$;
these exist since $w,w_s\in L^2$. For the other endpoints take
$j\to\infty$ at fixed $J$. The $C^1$ collar traces converge to
those of the three identified profiles, so the resulting boundary
terms vanish exponentially as $J\to\infty$. The tension integral on
the right of \eqref{eq:edge-coercivity} tends to zero even over the
whole half-line by \eqref{eq:strong-slices}. Thus the collars exhaust
the selected energy. Their limiting energies, as $J\to\infty$, are
$f_\pi E_b$, $4\pi$ and $4\pi$, respectively. Hence
$E(u(t_j))\to 8\pi+f_\pi E_b$; energy monotonicity gives the same
limit \eqref{eq:energy-limit} for all $t\to\infty$.

\paragraph{Convergence without selecting times.}

Set
\[
 A=\alpha_s+\sin\alpha,\qquad B=\beta_s+|\sin\beta|,\qquad
 H_+(\beta)=
 \begin{cases}1+\cos\beta,&\beta>\pi,\\0,&\beta\le\pi.\end{cases}
\]
The absolute value in $B$ accounts both for the decreasing overshoot
from $b$ to $\pi$ in the root and for the decreasing bubble from
$\pi$ to $0$. For the two unweighted scalar energies, completion of
squares gives the boundary contributions $4\pi$ and $4\pi+E_b$.
For the extra weight $\mathfrak f-1$ complete the square only on
$\{\beta>\pi\}$, using
$(H_+(\beta))_s=\mathbf 1_{\{\beta>\pi\}}\beta_s|\sin\beta|$.
Integration by parts gives an additional boundary contribution
$(f_\pi-1)E_b$. Altogether,
\begin{align}
 E-(8\pi+f_\pi E_b)
 ={}&\pi\|A\|_2^2+\pi\|B\|_2^2
 +\pi\int_{\beta>\pi}(\mathfrak f-1)B^2\,ds\notag\\
 &+\pi\int_{\beta\le\pi}(\mathfrak f-1)
                    (\beta_s^2+\sin^2\beta)\,ds\notag\\
 &+2\pi\int\mathfrak f'(\alpha)\alpha_s H_+(\beta)\,ds.
 \label{eq:sharp-squares}
\end{align}

Only the last term lacks a sign. In the fixed tube its integrand is
supported in $I_\alpha=[S_\alpha-K,S_\alpha+K]$ for a fixed $K$.
Let $v=\tan((\beta-\pi)_+/2)$. On the overshoot set,
$v_s+v=(1+v^2)B/2$; elsewhere $v_s+v=0$ almost everywhere.
Integration from $v(0)=\tan(\delta/2)$ and Cauchy--Schwarz give
\[
 \sup_{I_\alpha}(\beta-\pi)_+
 \le C\delta e^{-S_\alpha}+C\|B\|_2.
\]
Since $H_+(\beta)\le(\beta-\pi)_+^2/2$ and
$\int_{I_\alpha}|\alpha_s|\le C$, the last term of
\eqref{eq:sharp-squares} is bounded in absolute value by
$C\epsilon_w(\delta^2e^{-2S_\alpha}+\|B\|_2^2)$.
To control the first summand, $z=\pi-\alpha$ obeys
$z(0)=0$ and $z_s=\sin z-A\le z+|A|$. The tube gives
$z(S_\alpha)\ge c>0$, so integration and Cauchy--Schwarz yield
$c\le C e^{S_\alpha}\|A\|_2$. Taking $\epsilon_w$ and then
$\delta$ sufficiently small absorbs the unsigned term and proves
\begin{equation}
 E-(8\pi+f_\pi E_b)
 \ge c(\|A\|_2^2+\|B\|_2^2)
 +\pi\int_{\beta\le\pi}(\mathfrak f-1)
                       (\beta_s^2+\sin^2\beta)\,ds.
 \label{eq:sharp-floor}
\end{equation}
Every term on the right tends to zero as $t\to\infty$ by the full-time
energy limit. In particular $e^{-S_\alpha}\le C\|A\|_2\to0$; thus
both centers escape the disk boundary.

For clarity, consider any sequence $t\to\infty$. In a fixed or
translated chart the energy gives local uniform compactness. The
first-order equations $\alpha_s=-\sin\alpha+A$ and
$\beta_s=-|\sin\beta|+B$ show that a local uniform limit solves
the corresponding exact equations and that convergence of the
derivatives is strong in $L^2_{\rm loc}$. The fixed chart, with
data $\pi,b$ at $s=0$, therefore gives $(\pi,H_b(e^{-s}))$.
At its own center each angle gives a nonconstant whole-line solution
in the tube about $Q$. It takes an interior value in $(0,\pi)$
at $x=0$, so uniqueness for the first-order equation makes it a
translate of $Q$. The phase condition fixes the translation:
$\int(Q(x-h)-Q(x)) Z(x)\,dx$ has strict sign for $h\ne0$.

If $L$ stayed bounded on a subsequence with $L\to\ell$, the two
profiles in the outer chart would converge to $Q(x)$ and $Q(x-\ell)$.
On a compact interval where $\mathfrak f(Q(x))-1>0$, the latter
lies strictly between $0$ and $\pi$. The interaction integral in
\eqref{eq:sharp-floor} would then have a positive lower limit,
contrary to that bound. Thus $L\to\infty$. In the outer chart the
second angle is consequently confined to a small neighborhood of
$\pi$ on the whole limiting line; in the inner chart the first
angle is similarly confined near $0$. A whole-line solution of
$y'=-|\sin y|$ confined near $\pi$ is constant: it is monotone
and both end limits must be zeros of the right-hand side.
The same argument applies to $y'=-\sin y$ in $[0,\pi/4]$.
This proves convergence of \emph{both} components in all three
charts. Since the initial sequence of times was arbitrary and the
limits are unique, the local strong energy and uniform convergence
hold in full time.

Finally, for fixed $J$ the energy on the root collar $[0,J]$ and
the two vertex collars $[S_i-J,S_i+J]$ converges to the energy of
the corresponding limiting profiles there. These collars are
eventually disjoint. Subtracting their energies from
\eqref{eq:energy-limit} and then sending $J\to\infty$ proves the
asserted double-limit vanishing of the energy on the three
complementary intervals. On each such interval the squared
$H^1$ norms of \emph{both} signed pole coordinates are bounded by
$C E(u;I_J^k\times S^1)$, since $\sin^2 w\ge c w^2$ in the fixed
pole neighborhoods. For fixed $J$, the finite intervals eventually
have length at least one; the $H^1$-to-$L^\infty$ inequality there
and on the half-line has a common constant (apply it on unit
subintervals). At any angular parameter $\theta$, a pair with signed
pole coordinates $w_\alpha,w_\beta$ lies at target distance at most
$|w_\alpha|+\sqrt{f_\pi}|w_\beta|$ from the pole pair, by moving
along the two meridians. Its entire angular orbit satisfies the
same estimate. The triangle inequality therefore bounds each full
image diameter by $C E(u;I_J^k\times S^1)^{1/2}$. This proves the
full-time strong no-neck statement of Theorem~\ref{thm:example}.

\subsection{The middle neck and the compactified domain}\label{sec:geometry}

Return to the selected times of \eqref{eq:strong-slices}. We now know
that $L(t_j)\to\infty$. By \eqref{eq:residue-tension}, $L\le Y$, and
$q_m\le q_\beta$,
\[
 |d(t_j)|L(t_j)
 \le C\sqrt{Y(t_j)^2\Diss(t_j)/q_m(t_j)}
                   +CL(t_j)e^{-cL(t_j)}\longrightarrow0 .
\]
Apply \eqref{eq:trimmed-packet-estimate} first with $K>1$ fixed and
$j\to\infty$; the heat loss multiplied by $L$ vanishes by
$\lambda\le e^{-L}$. Then let $K\to\infty$. The two remaining errors,
$CK^{-2}$ and $Ce^{-cK^2}$, both vanish. This proves the last
assertion of Theorem~\ref{thm:example}. In fact, in $s=-\log r$,
\[
 q(t_j,r_m)\,e^{2\phi_u}(dr^2+r^2d\theta^2)
 =\frac{q(t_j,r_m)}{q(t_j,r)}(ds^2+d\theta^2),
\]
so the comparison is with the unit cylinder after \emph{one}
normalization on each selected slice. The two packet residues need not
be equal. This argument does not imply full-time convergence of the metric.

\paragraph{The compactified domain.}
This last observation is topological. Use annular charts with their nonterminal ends represented by boundary
circles, and keep the original outer boundary of the disk. In logarithmic
polar coordinates, the charts are products of intervals with $S^1$. Gluing
in radial order along the circles concatenates these intervals and leaves
a punctured disk. Adding one point
at the terminal end gives a space homeomorphic to the closed disk,
with its original outer boundary. Collapsing the two interface circles
instead gives the nodal picture with two spheres. Only the middle
annulus is the subject of the metric limit above.

\end{document}